\documentclass[11pt]{article}
\usepackage{mathrsfs}
\usepackage{amsfonts}
\usepackage{amssymb}
\usepackage{amsmath}
\usepackage[all]{xy}

\def\Hom{\mathop{\rm Hom}\nolimits}
\def\Tor{\mathop{\rm Tor}\nolimits}
\def\Coker{\mathop{\rm Coker}\nolimits}
\def\Mod{\mathop{\rm Mod}\nolimits}
\def\mod{\mathop{\rm mod}\nolimits}
\def\pd{\mathop{\rm pd}\nolimits}
\def\id{\mathop{\rm id}\nolimits}
\def\fd{\mathop{\rm fd}\nolimits}
\def\Im{\mathop{\rm Im}\nolimits}
\def\Tr{\mathop{\rm Tr}\nolimits}
\def\Ext{\mathop{\rm Ext}\nolimits}

\def\Gdim{\mathop{\rm G-dim}\nolimits}
\def\dim{\mathop{\rm dim}\nolimits}

\def\mod{\mathop{\rm mod}\nolimits}
\title{\Large \bf Torsionfree Dimension of Modules and Self-Injective
Dimension of Rings
\thanks{This research was partially supported by the
Specialized Research Fund for the Doctoral Program of Higher
Education (Grant No. 20100091110034), NSFC (Grant No. 10771095) and NSF
of Jiangsu Province of China (Grant Nos. BK2010047, BK2010007).}}
\author{Chonghui Huang and Zhaoyong Huang}
\date{}
\begin{document}
\baselineskip=18pt
\maketitle

\begin{abstract}
Let $R$ be a left and right Noetherian ring. We introduce the notion
of the torsionfree dimension of finitely generated $R$-modules. For
any $n\geq 0$, we prove that $R$ is a Gorenstein ring with
self-injective dimension at most $n$ if and only if every finitely
generated left $R$-module and every finitely generated right
$R$-module have torsionfree dimension at most $n$, if and only if
every finitely generated left (or right) $R$-module has Gorenstein
dimension at most $n$. For any $n \geq 1$, we study the properties
of the finitely generated $R$-modules $M$ with $\Ext _R^i(M, R)=0$
for any $1\leq i \leq n$. Then we investigate the relation between
these properties and the self-injective dimension of $R$.
\end{abstract}

\vspace{0.5cm}

\centerline{\bf 1. Introduction}

\vspace{0.2cm}

Throughout this paper, $R$ is a left and right Noetherian ring
(unless stated otherwise) and $\mod R$ is the category of finitely
generated left $R$-modules. For a module $M\in \mod R$, we use $\pd
_R M$, $\fd _R M$, $\id _R M$ to denote the projective, flat,
injective dimension of $M$, respectively.

For any $n\geq 1$, we denote $^{\bot_n}{_RR}=\{M\in \mod R\ |\ \Ext
_R^i(_RM, {_RR})=0$ for any $1\leq i \leq n\}$ (resp.
$^{\bot_n}R_R=\{N\in \mod R^{op}\ |\ \Ext^i_{R^{op}}(N_R, {R_R})=0$
for any $1\leq i \leq n\}$), and $^{\bot}{_RR}=\bigcap _{n\geq
1}{^{\bot_n}{_RR}}$ (resp. $^{\bot}R_R=\bigcap _{n\geq
1}{^{\bot_n}R_R}$).

For any $M\in \mod R$, there exists an exact sequence: $$P_{1}
\buildrel {f} \over \longrightarrow P_{0} \to M \to 0$$ in $\mod R$
with $P_{0}$ and $P_{1}$ projective. Then we get an exact sequence:
$$0 \to M^* \to P_{0}^* \buildrel {f^*} \over \longrightarrow
P_{1}^* \to \Tr M \to 0$$ in $\mod R^{op}$, where $(-)^*= \Hom(-,
R)$ and $\Tr M=\Coker f^*$ is the {\it transpose} of $M$ ([1]).

Auslander and Bridger generalized the notions of finitely generated
projective modules and the projective dimension of finitely
generated modules as follows.

\vspace{0.2cm}

{\bf Definition 1.1} ([1]). Let $M \in \mod R$.

(1) $M$ is said to have {\it Gorenstein dimension zero} if $M\in
{^{\bot}{_RR}}$ and $\Tr M \in {^{\bot}R_R}$.

(2) For a non-negative integer $n$, the {\it Gorenstein dimension}
of $M$, denoted by $\Gdim _RM$, is defined as inf$\{ n\geq 0\ |\
$there exists an exact sequence $0 \to M_{n} \to \cdots \to M_{1}
\to M_{0} \to M \to 0$ in $\mod R$ with all $M_i$ having Gorenstein
dimension zero$\}$. We set $\Gdim _RM$ infinity if no such integer
exists.

\vspace{0.2cm}

Huang introduced in [7] the notion of the left orthogonal
dimension of modules as follows, which is ``simpler" than that of
the Gorenstein dimension of modules.

\vspace{0.2cm}

{\bf Definition 1.2} ([7]). For a module $M\in \mod R$, the {\it
left orthogonal dimension} of a module $M \in \mod R$, denoted by
${^{\bot}{_RR}}-\dim _RM$, is defined as inf$\{ n\geq 0\ |\ $there
exists an exact sequence $0 \to X_{n} \to \cdots \to X_{1} \to X_{0}
\to M \to 0$ in $\mod R$ with all $X_i\in {^{\bot}{_RR}}\}$. We set
${^{\bot}{_RR}}-\dim _RM$ infinity if no such integer exists.

\vspace{0.2cm}

Let $M\in \mod R$. It is trivial that ${^{\bot}{_RR}}-\dim _RM \leq
\Gdim _RM$. On the other hand, by [14], we have that
${^{\bot}{_RR}}-\dim _RM\neq \Gdim _RM$ in general.

Recall that $R$ is called a {\it Gorenstein ring} if $\id _RR=\id
_{R^{op}}R<\infty$. The following result was proved by Auslander and
Bridger in [1, Theorem 4.20] when $R$ is a commutative Noetherian
local ring. Hoshino developed in [4] Auslander and Bridger's
arguments and applied obtained the obtained results to Artinian
algebras. Then Huang generalized in [7, Corollary 3] Hoshino's
result with the left orthogonal dimension replacing the Gorenstein
dimension of modules.

\vspace{0.2cm}

{\bf Theorem 1.3} ([4, Theorem] and [7, Corollary 3]). {\it The
following statements are equivalent for an Artinian algebra $R$.

(1) $R$ is Gorenstein.

(2) Every module in $\mod R$ has finite Gorenstein dimension.

(3) Every module in $\mod R$ and every module in $\mod R^{op}$ have
finite left orthogonal dimension.}

\vspace{0.2cm}

One aim of this paper is to generalize this result to left and right
Noetherian rings. On the other hand, note that the left orthogonal
dimension of modules is defined by the least length of the
resolution composed of the modules in ${^{\bot}{_RR}}$, which are
the modules satisfying one of the two conditions in the definition
of modules having Gorenstein dimension zero. So, a natural question
is: If a new dimension of modules is defined by the least length of
the resolution composed of the modules satisfying the other
condition in the definition of modules having Gorenstein dimension
zero, then can one give an equivalent characterization of Gorenstein
rings similar to the above result in terms of the new dimension of
modules? The other aim of this paper is to give a positive answer to
this question. This paper is organized as follows.

In Section 2, we give the definition of $n$-torsionfree modules, and
investigate the properties of such modules. We prove that a module
in $\mod R$ is $n$-torsionfree if and only if it is an $n$-syzygy of
a module in $^{\bot_n}{_RR}$.

In Section 3, we introduce the notion of the torsionfree dimension
of modules. Then we give some equivalent characterizations of
Gorenstein rings in terms of the properties of the torsionfree
dimension of modules. The following is the main result in this
paper.

\vspace{0.2cm}

{\bf Theorem 1.4.} {\it For any $n\geq 0$, the following statements
are equivalent.

(1) $R$ is a Gorenstein ring with $\id _RR=\id _{R^{op}}R\leq n$.

(2) Every module in $\mod R$ has Gorenstein dimension at most $n$.

(3) Every module in $\mod R^{op}$ has Gorenstein dimension at most
$n$.

(4) Every module in $\mod R$ and every module in $\mod R^{op}$ have
torsionfree dimension at most $n$.

(5) Every module in $\mod R$ and every module in $\mod R^{op}$ have
left orthogonal dimension at most $n$.}

\vspace{0.2cm}

In Section 4, for any $n\geq 1$, we first prove that every module in
${^{\bot_n}{_RR}}$ is torsionless (in this case, ${^{\bot_n}{_RR}}$
is said to have the {\it torsionless property}) if and only if every
module in ${^{\bot_n}{_RR}}$ is $\infty$-torsionfree, if and only if
every module in ${^{\bot _n}{_RR}}$ has torsionfree dimension at
most $n$, if and only if every $n$-torsionfree module in $\mod R$ is
$\infty$-torsionfree, if and only if every $n$-torsionfree module in
$\mod R^{op}$ is in $ {^{\bot}R_R}$, if and only if
${^{\bot_n}R_R}={^{\bot}R_R}$. Note that if $\id _{R^{op}}R\leq n$,
then ${^{\bot_n}{_RR}}$ has the torsionless property. As some
applications of the obtained results, we investigate when the
converse of this assertion holds true. Assume that $n$ and $k$ are
positive integers and $^{\bot_n}{_RR}$ has the torsionless property.
If $R$ is $g_n(k)$ or $g_n(k)^{op}$ (see Section 4 for the
definitions), then $\id _{R^{op}} R \leq n+k-1$. As a corollary, we
have that if $\id _R R\leq n$, then $\id _R R=\id _{R^{op}}R\leq n$
if and only if $^{\bot _n}{_RR}$ has the torsionless property.

In view of the results obtained in this paper, we pose in Section 5
the following two questions: (1) Is the subcategory of $\mod R$
consisting of modules with torsionfree dimension at most $n$ closed
under extensions or under kernels of epimorphisms? (2) If $\id
_{R^{op}} R\leq n$, does then every module $M\in \mod R$ has
torsionfree dimension at most $n$?

\vspace{0.5cm}

\centerline{\bf 2. Preliminaries}

\vspace{0.2cm}

Let $M\in \mod R$ and $n\geq 1$. Recall from [1] that $M$ is called
{\it $n$-torsionfree} if $\Tr M \in {^{\bot _n}R_R}$; and $M$ is
called {\it $\infty$-torsionfree} if $M$ is $n$-torsionfree for all
$n$. We use $\mathcal{T} _n (\mod R)$ (resp. $\mathcal{T} (\mod R)$)
to denote the subcategory of $\mod R$ consisting of all
$n$-torsionfree modules (resp. $\infty$-torsionfree modules). It is
well-known that $M$ is 1-torsionfree (resp. 2-torsionfree) if and
only if $M$ is torsionless (resp. reflexive) (see [1]). Also recall
from [1] that $M$ is called an $n$-syzygy module (of $A$), denoted
by $\Omega^n(A)$, if there exists an exact sequence $0\to M\to
P_{n-1}\to \cdots \to P_1\to P_0\to A\to 0$ in $\mod R$ with all
$P_i$ projective. In particular, set $\Omega^0(M)=M$. We use $\Omega
^n (\mod R)$ to denote the subcategory of $\mod R$ consisting of all
$n$-syzygy modules. It is easy to see that $\mathcal {T}_n (\mod R)
\subseteq \Omega ^n (\mod R)$, and in general, this inclusion is
strict when $n\geq 2$ (see [1]).

Jans proved in [13, Corollary 1.3] that a module in $\mod R$ is
1-torsionfree if and only if it is an 1-syzygy of a module in
$^{\bot_1}{_RR}$. We generalize this result as follows.

\vspace{0.2cm}

{\bf Proposition 2.1.} {\it For any $n \geq 1$, a module in $\mod R$
is $n$-torsionfree if and only if it is an $n$-syzygy of a module in
$^{\bot_n}{_RR}$}.

\vspace{0.2cm}

{\bf Proof}. Assume that $M\in \mod R$ is an $n$-syzygy of a module
$A$ in $^{\bot_n}{_RR}$. Then there exists an exact sequence:
$$0\to M\to P_{n-1}\to\cdots
 \to P_1 \buildrel {f} \over \longrightarrow P_0
\to A \to 0$$ in $\mod R$ with all $P_i$ projective. Let $$P_{n+1}
\to P_n \to  M \to 0$$ be a projective presentation of $M$ in $\mod
R$. Then the above two exact sequences yield the following exact
sequence:
$$0\to A^*\to P^*_0\buildrel {f^*} \over \longrightarrow\cdots
 \to P^*_n  \to P^*_{n+1}
\to \Tr M \to 0.$$ By the exactness of $P_{n+1}\to P_n\to\cdots
 \to P_1 \buildrel {f} \over \longrightarrow P_0 $, we
get that $\Tr M\in {^{\bot_n}R_R}$. Thus $M$ is $n$-torsionfree.

Conversely, assume that $M\in \mod R$ is $n$-torsionfree and $$P_1
\buildrel {g} \over\longrightarrow P_0\buildrel {\pi}
\over\longrightarrow  M \to 0$$ is a projective presentation of
$M\in \mod R$. Then we get an exact sequence: $$0\to M^* \buildrel
{\pi ^*} \over\longrightarrow P^*_0\buildrel {g^*} \over
\longrightarrow P^*_1 \to \Tr M \to 0$$ in $\mod R^{op}$. Let
$$\cdots \buildrel {h_{n+1}} \over\longrightarrow Q_{n}
\buildrel {h_n} \over\longrightarrow \cdots\buildrel {h_1}
\over\longrightarrow Q_{0}\buildrel {h_0} \over\longrightarrow M^*
\to 0$$ be a projective resolution of $M^*$ in $\mod R^{op}$. Then
we get a projective resolution of $\Tr M$:
$$\cdots \buildrel {h_{n+1}} \over\longrightarrow Q_{n}
\buildrel {h_n} \over\longrightarrow \cdots\buildrel {h_1}
\over\longrightarrow Q_{0} \buildrel {\pi^* h_0}
\over\longrightarrow P^*_0 \buildrel {g^*} \over\longrightarrow
P^*_1 \to \Tr M\to 0.$$ Because $M$ is $n$-torsionfree, $\Tr M\in
{^{\bot_n}R_R}$ and we get the following exact sequence:
$$0\to (\Tr M)^*\to P^{**}_1\buildrel {g^{**}} \over\longrightarrow
P^{**}_0\buildrel {h_0^*\pi^{**}} \over\longrightarrow
Q^*_0\buildrel {h_1^{*}} \over\longrightarrow\cdots
 \buildrel {h_{n-1}^{*}} \over\longrightarrow Q^*_{n-1}
\to \Coker h_{n-1}^{*} \to 0.$$ It is easy to see that $M\cong
\Coker g^{**}$. By the exactness of $Q_{n-1} \buildrel {h_{n-1}}
\over\longrightarrow \cdots\buildrel {h_1} \over\longrightarrow
Q_{0} \buildrel {\pi^* h_0} \over\longrightarrow P^*_0\buildrel
{g^*} \over \longrightarrow P^*_1$, we get that $\Coker h_{n-1}^*\in
{^{\bot_n}{_RR}}$. The proof is finished. \hfill $\square$

\vspace{0.2cm}

As an immediate consequence, we have the following

\vspace{0.2cm}

{\bf Corollary 2.2.} {\it For any $n \geq 1$, an $n$-torsionfree
module in $\mod R$ is a 1-syzygy of an $(n-1)$-torsionfree module
$A$ in $\mod R$ with $A \in {^{\bot_1}{_RR}}$. In particular, an
$\infty$-torsionfree module in $\mod R$ is a 1-syzygy of an
$\infty$-torsionfree module $T$ in $\mod R$ with $T \in
{^{\bot_1}{_RR}}$.}

\vspace{0.2cm}

We also need the following easy observation.

\vspace{0.2cm}

{\bf Lemma 2.3.} {\it For any $n \geq 1$, both $\mathcal{T}_n(\mod
R)$ and $\mathcal{T}(\mod R)$ are closed under direct summands and
finite direct sums.}

\vspace{0.5cm}

\centerline{\bf 3. Torsionfree dimension of modules}

\vspace{0.2cm}

In this section, we will introduce the notion of the torsionfree
dimension of modules in $\mod R$. Then we will give some equivalent
characterizations of Gorenstein rings in terms of the properties of
this dimension of modules.

We begin with the following well-known observation.

\vspace{0.2cm}

{\bf Lemma 3.1} ([1, Lemma 3.9]). {\it Let $0\to A\buildrel {f}
\over\longrightarrow B\to C\to 0$ be an exact sequence in $\mod R$.
Then we have exact sequences $0\to C^*\to B^*\to A^*\to \Coker
f^*\to 0$ and $0\to \Coker f^*\to \Tr C\to \Tr B\to \Tr A\to 0$ in
$\mod R^{op}$.}

\vspace{0.2cm}

The following result is useful in this section.

\vspace{0.2cm}

{\bf Proposition 3.2.} {\it Let $$0\to M\to T_1\buildrel {f}
\over\longrightarrow T_0\to A\to 0$$ be an exact sequence in $\mod
R$ with both $T_0$ and $T_1$ in $\mathcal{T}(\mod R)$. Then there
exists an exact sequence:
$$0\to M\to P \to T\to A\to 0$$ in $\mod R$ with $P$ projective and
$T\in \mathcal{T}(\mod R)$.}

\vspace{0.2cm}

{\it Proof.} Let $$0\to M\to T_1\buildrel {f} \over\longrightarrow
T_0\to A\to 0$$ be an exact sequence in $\mod R$ with both $T_0$ and
$T_1$ in $\mathcal{T}(\mod R)$. By Corollary 2.2, there exists an
exact sequence $0\to T_1 \to P \to W\to 0$ in $\mod R$ with $P$
projective and $W\in {^{\bot _1}{_RR}}\bigcap \mathcal{T}(\mod R)$.
Then we have the following push-out diagram:
$$\xymatrix{
& & 0 \ar[d] &0 \ar[d] &   \\
0 \ar[r]  & M \ar@{=}[d]  \ar[r] & T_1 \ar[d]
\ar[r]  &\Im f \ar[d] \ar[r] & 0 \\
0 \ar[r] & M \ar[r] & P \ar[d] \ar[r] & B \ar[d] \ar[r] & 0  \\
&  & W \ar[d]
\ar@{=}[r] & W\ar[d]  &   \\
&  & 0  & 0   &   } $$ Now, consider the following push-out diagram:
$$\xymatrix{
& 0 \ar[d] &0 \ar[d] &   &    \\
0 \ar[r] & \Im f  \ar[d] \ar[r] &T_0
\ar[d]  \ar[r]  & A \ar@{=}[d] \ar[r]  &0  \\
0 \ar[r]  & B \ar[d]  \ar[r] &T
\ar[d] \ar[r] & A \ar[r]  &0  \\
& W\ar[d]\ar@{=}[r] &W \ar[d] &   &    \\
&  0 & 0  &  & }$$ Because $W\in {^{\bot _1}{_RR}}$, we get an exact
sequence:
$$0\to \Tr W \to \Tr T \to \Tr T_0 \to 0$$ by Lemma 3.1 and the
exactness of the middle column in the above diagram. Because both
$W$ and $T_0$ are in $\mathcal{T}(\mod R)$, both $\Tr W$ and $\Tr
T_0$ are in ${^{\bot}R_R}$. So $\Tr T$ is also in ${^{\bot}R_R}$ and
hence $T\in\mathcal{T}(\mod R)$. Connecting the middle rows in the
above two diagrams, then we get the desired exact sequence. \hfill
$\square$

\vspace{0.2cm}

Now we introduce the notion of the torsionfree dimension of modules
as follows.

\vspace{0.2cm}

{\bf Definition 3.3.} For a module $M\in \mod R$, the {\it
torsionfree dimension} of $M$, denoted by $\mathcal{T}-\dim _RM$, is
defined as inf$\{ n\geq 0\ |\ $there exists an exact sequence $0 \to
X_{n} \to \cdots \to X_{1} \to X_{0} \to M \to 0$ in $\mod R$ with
all $X_i\in \mathcal{T}(\mod R)\}$. We set $\mathcal{T}-\dim _RM$
infinity if no such integer exists.

\vspace{0.2cm}

Let $M\in \mod R$. It is trivial that $\mathcal{T}-\dim _RM \leq
\Gdim _RM$. On the other hand, by [14], we have that
$\mathcal{T}-\dim _RM\neq \Gdim _RM$ in general.


\vspace{0.2cm}

{\bf Proposition 3.4.} {\it Let $M\in \mod R$ and $n\geq 0$. If
$\mathcal {T}-\dim _R M\leq n$, then there exists an exact sequence
$0\to H\to T\to M\to 0$ in $\mod R$ with $\pd _R H\leq n-1$ and
$T\in\mathcal{T}(\mod R)$.}

\vspace{0.2cm}

{\bf Proof.} We proceed by induction on $n$. If $n=0$, then $H=0$
and $T=M$ give the desired exact sequence. If $n=1$, then there
exists an exact sequence:
$$0\to T_1\to T_0\to M\to 0$$ in $\mod R$ with both $T_0$ and $T_1$
in $T\in\mathcal{T}(\mod R)$. Applying Proposition 3.2, with $A=0$,
we get an exact sequence:
$$0\to P\to T'_0\to M\to 0 $$ in $\mod R$ with $P$ projective and
$T'_0\in \mathcal {T}(\mod R)$.

Now suppose $n\geq 2$. Then there exists an exact sequence: $$0\to
T_n\to T_{n-1}\to\cdots \to T_0\to M\to 0$$ in $\mod R$ with all
$T_i\in \mathcal{T}(\mod R)$. Set $K=\Im(T_{1}\to T_{0})$. By the
induction hypothesis, we get the following exact sequence:
$$0\to P_n\to P_{n-1}\to P_{n-2}\to \cdots \to P_{2}\to T'_1\to K\to 0$$
in $\mod R$ with all $P_i$ projective and $T^{'}_1\in
\mathcal{T}(\mod R)$. Set $N=\Im(P_2\to T'_1)$. By Proposition 3.2,
we get an exact sequence: $$0\to N\to P_1\to T\to M\to 0$$ in $\mod
R$ with $P_1$ projective and $T\in \mathcal{T}(\mod R)$. Thus we get
the desired exact sequence: $$0\to P_n\to P_{n-1}\to P_{n-2}\to
\cdots \to P_{1}\to T\to M\to 0$$ and the assertion follows. \hfill
$\square$

\vspace{0.2cm}

Christensen, Frankild and Holm proved in [2, Lemma 2.17] that a
module with Gorenstein dimension at most $n$ can be embedded into a
module with projective dimension at most $n$, such that the cokernel
is a module with Gorenstein dimension zero. The following result
extends this result.

\vspace{0.2cm}

{\bf Corollary 3.5.} {\it Let $M\in \mod R$ and $n\geq 0$. If
$\mathcal {T}-\dim _R M\leq n$, then there exists an exact sequences
$0\to M\to N\to T\to 0$ in $\mod R$ with $\pd _R N\leq n$ and $T\in
{^{\bot _1}{_RR}} \bigcap \mathcal{T}(\mod R)$.}

\vspace{0.2cm}

{\it Proof.} Let $M\in \mod R$ with $\mathcal {T}-\dim _R M\leq n$.
By Proposition 3.4, there exists an exact sequence $0\to H\to
T^{'}\to M\to 0$ in $\mod R$ with $\pd _R H\leq n-1$ and
$T^{'}\in\mathcal{T}(\mod R)$. By Corollary 2.2, there exists an
exact sequence $0\to T^{'}\to P\to T\to 0$ in $\mod R$ with $P$
projective and $T\in {^{\bot _1}{_RR}} \bigcap \mathcal{T}(\mod R)$.
Consider the following push-out diagram:
$$\xymatrix{
& & 0 \ar[d] &0 \ar[d] &   \\
0 \ar[r]  & H \ar@{=}[d]  \ar[r] & T^{'} \ar[d]
\ar[r]  &M \ar[d] \ar[r] & 0 \\
0 \ar[r] & H \ar[r] & P \ar[d] \ar[r] & N \ar[d] \ar[r] & 0  \\
&  & T \ar[d]
\ar@{=}[r] & T\ar[d]  &   \\
&  & 0  & 0   &   } $$ Then the third column in the above diagram is
as desired.  \hfill $\square$

\vspace{0.2cm}

The following result plays a crucial role in proving the main result
in this paper.

\vspace{0.2cm}

{\bf Theorem 3.6.} {\it For any $n\geq 0$, if every module in $\mod
R$ has torsionfree dimension at most $n$, then $\id _{R^{op}} R\leq
n$.}

\vspace{0.2cm}

To prove this theorem, we need some lemmas. We use $\Mod R$ to
denote the category of left $R$-modules.

\vspace{0.2cm}

{\bf Lemma 3.7} ([11, Proposition 1]). {\it $\id _{R^{op}} R=\sup\{\fd
_RE\ |\ E$ is an injective module in $\Mod R\}=\fd _RQ$ for any
injective cogenerator $Q$ for $\Mod R$.}

\vspace{0.2cm}

{\bf Lemma 3.8.} {\it For any $n \geq 0$, $\id _{R^{op}} R\leq n$ if
and only if every module in $\mod R$ can be embedded into a module
in $\Mod R$ with flat dimension at most $n$.}

\vspace{0.2cm}

{\it Proof.} Assume that $\id _{R^{op}} R\leq n$. Then the injective
envelope of any module in $\mod R$ has flat dimension at most $n$ by
Lemma 3.7, and the necessity follows.

Conversely, let $E$ be any injective module in $\Mod R$. Then by
[15, Exercise 2.32], $E=\underset{i\in I}{\underset{\rightarrow}{\rm
lim}}M_i$, where $\{M_i\ |\ i \in I\}$ is the set of all finitely
generated submodules of $E$ and $I$ is a directed index set (in
which the quasi-order is defined by $i\leq j$ if and only if $M_i
\leq M_j$, the homomorphism $\lambda^i_j:\ M_i \to M_j$ is the
canonical embedding). By assumption, for any $i\in I$ and $M_i\in
\mod R$, we have an exact sequence $0\to M_i\buildrel {\alpha _i}
\over\to N_i$ with $N_i\in\Mod R$ and $\fd _RN_i\leq n$.

Put $K=\underset{i\in I}{{\prod}}N_i$ and $I_i$=$\{j\in I\ |\
M_i\leq M_j \}$ for any $i\in I$. Since $R$ is a left and right
Noetherian ring, any direct product of flat modules is still flat. So $\fd
_RK\leq n$. Define $\beta_i=\underset{k\in I}{{\prod}}f^i_k$ with
$$f^i_k=\left\{\aligned
&  \alpha_k\lambda^i_k,\ {\rm if} \ k\in I_i,\\
& 0, \ \ \ \ \ \ \ {\rm if} \ k\notin I_i \endaligned\right. $$
for any $i, k\in I$. Then $0\to M_i \buildrel {\beta _i} \over\to K$ is
exact for any $i\in I$. For any $i\leq j$ (determined
by $M_i\leq M_j$), we have the following commutative and exact
diagram:
$$\xymatrix{
& 0 \ar[d] & &  \\
0 \ar[r] & M_i  \ar[d]^{\lambda ^i_j} \ar[r]^{\beta_i} &K
\ar[d]^{\varphi_j^i}  &    \\
0 \ar[r]  & M_j \ar[r]^{\beta_j} &K  &  \\
 }$$ where $\varphi_j^i=\underset{k\in I}{{\prod}}h_k$ with
$$h_k=\left\{\aligned
&  1_{N_k},\ {\rm if} \ k\in I_j,\\
& 0, \ \ \ \ {\rm if} \ k\notin I_j \endaligned\right.$$ for any $k\in I$. It is
clear that $\{K, \varphi_j^i\}$ is a direct system of the constant module $K$.
It follows from [15, Theorem 2.18] that we get a monomorphism
$0\to E(=\underset{i\in I}{\underset{\rightarrow}{\rm lim}}M_i)\to
\underset{i\in I}{\underset{\rightarrow}{\rm lim}}K$. Because the
functor $\Tor$ commutes with $\underset{i\in
I}{\underset{\rightarrow}{\rm lim}}$ by [15, Theorem 8.11], $\fd
_R\underset{i\in I}{\underset{\rightarrow}{\rm lim}}K \leq n$. So
$\fd _RE\leq n$ and hence $\id _{R^{op}} R\leq n$ by Lemma 3.7.
\hfill $\square$

\vspace{0.2cm}

{\it Proof of Theorem 3.6.} By assumption and Corollary 3.5, we have
that every module in $\mod R$ can be embedded into a module in $\mod
R$ with projective dimension at most $n$. Then by Lemma 3.8, we get
the assertion. \hfill $\square$

\vspace{0.2cm}

{\bf Lemma 3.9.} {\it For any $M\in \mod R$ and $n\geq 0$,
${^{\bot}{_RR}}-\dim _RM\leq n$ if and only if $\Ext _R^{n+i}(M,
R)=0$ for any $i \geq 1$.}

\vspace{0.2cm}

{\it Proof.} For any $M\in \mod R$, consider the following exact
sequence:
$$\cdots \to W_n \to W_{n-1} \to \cdots \to W_0 \to M \to 0$$
in $\mod R$ with all $W_i$ in ${^{\bot}{_RR}}$. Then we have that
$\Ext _R^{i}(\Im (W_n \to W_{n-1}), R)\cong \Ext _R^{n+i}(M, R)$ for
any $i\geq 1$. So $\Im (W_n \to W_{n-1}) \in {^{\bot}{_RR}}$ if and
only if $\Ext _R^{n+i}(M, R)=0$ for any $i \geq 1$, and hence the
assertion follows. \hfill $\square$

\vspace {0.2cm}

{\bf Proposition 3.10.} {\it For any $n\geq 0$, every module in $\mod
R$ has left orthogonal dimension at most $n$ if and only if $\id
_RR\leq n$.}

\vspace {0.2cm}

{\it Proof.} By Lemma 3.9, we have that $\id _RR\leq n$ if and only
if $\Ext _R^{n+i}(M, R)=0$ for any $M\in \mod R$ and $i\geq 1$, if
and only if ${^{\bot}{_RR}}-\dim _RM\leq n$ for any $M\in \mod R$.
\hfill $\square$

\vspace {0.2cm}

{\it Proof of Theorem 1.4.} $(1)\Rightarrow (2)+(3)$ follows from
[10, Theorem 3.5].

$(2)\Rightarrow (1)$ Let $M$ be any module in $\mod R$. Then by
assumption, we have that $\Gdim _RM$\linebreak $\leq n$ and
$\mathcal{T}-\dim _RM\leq n$. So $\id _{R^{op}}R\leq n$ by Theorem
3.6. On the other hand, because ${^{\bot}{_RR}}-\dim _RM\leq \Gdim
_RM$, $\id _RR\leq n$ by Proposition 3.10.

Symmetrically, we get $(3)\Rightarrow (1)$.

$(4)\Rightarrow (1)$ By Theorem 3.6 and its symmetric version.

$(2)+(3)\Rightarrow (4)$ Because $\mathcal{T}-\dim _RM\leq \Gdim
_RM$ and $\mathcal{T}-\dim _{R^{op}}N\leq \Gdim _{R^{op}}N$ for any
$M\in \mod R$ and $N\in \mod R^{op}$, the assertion follows.

$(1)\Leftrightarrow (5)$ By Proposition 3.10 and its symmetric
version. \hfill $\square$

\vspace{0.5cm}

\centerline{\bf 4. The torsionless property and self-injective
dimension}

\vspace{0.2cm}

The following result plays a crucial role in this section, which
generalizes [4, Lemma 4], [8, Lemma 2.1] and [13, Theorem 5.1].

\vspace{0.2cm}

{\bf Proposition 4.1.} {\it For any $n\geq 1$, the
following statements are equivalent.

(1) ${^{\bot_n}{_RR}}\subseteq \mathcal{T}_1(\mod R)$. In this case,
${^{\bot_n}{_RR}}$ is said to have the torsionless property.

(2) ${^{\bot_n}{_RR}}\subseteq \mathcal{T}(\mod R)$.

(3) Every module in ${^{\bot _n}{_RR}}$ has torsionfree dimension at
most $n$.

(4) $\mathcal{T}_n(\mod R)=\mathcal{T}(\mod R)$.

(5) $\mathcal{T}_n(\mod R^{op})\subseteq {^{\bot}R_R}$.

(6) ${^{\bot_n}R_R}={^{\bot}R_R}$.}

\vspace{0.2cm}

{\it Proof.} $(2)\Rightarrow (1)$ and $(2)\Rightarrow (3)$ are
trivial, and $(1)\Leftrightarrow (6)$ follows from [8, Lemma 2.1].
Note that $M$ and $\Tr\Tr M$ are projectively equivalent for any
$M\in \mod R$ or $\mod R^{op}$. Then it is not difficult to verify
$(2)\Leftrightarrow (5)$ and $(4)\Leftrightarrow (6)$. So it
suffices to prove $(1)\Rightarrow (2)$ and $(3)\Rightarrow (2)$.

$(1)\Rightarrow (2)$ Assume that $M\in {^{\bot_n}{_RR}}$. Then $M$
is torsionless by (1). So, by Proposition 2.1, we have an exact
sequence $0\to M \to P_0 \to M_1 \to 0$ in $\mod R$ with $P_0$
projective and $M_1 \in {^{\bot_1}{_RR}}$, which yields that $M_1\in
{^{\bot_{n+1}}{_RR}}$. Then $M_1$ is torsionless by (1), and again
by Proposition 2.1, we have an exact sequence $0\to M_1 \to P_1 \to
M_2 \to 0$ in $\mod R$ with $P_1$ projective and $M_2 \in
{^{\bot_1}{_RR}}$, which yields that $M_1\in {^{\bot_{n+2}}{_RR}}$.
Repeating this procedure, we get an exact sequence:
$$0\to M \to P_0 \to P_1 \to \cdots \to P_i \to \cdots$$
in $\mod R$ with all $P_i$ projective and $\Im (P_i \to P_{i+1})\in
{^{\bot_{n+i+1}}{_RR}}\subseteq {^{\bot_{i+1}}{_RR}}$, which implies
that $M$ is $\infty$-torsionfree by Proposition 2.1.

$(3)\Rightarrow (2)$ Assume that $M\in {^{\bot _n}{_RR}}$. Then
$\mathcal{T}-\dim _RM\leq n$ by assumption. By Proposition 3.4,
there exists an exact sequence:
$$0\to H\to T\to M\to 0\eqno{(1)}$$ in $\mod R$ with $\pd
_R H\leq n-1$ and $T\in\mathcal{T}(\mod R)$. Because $M\in {^{\bot
_n}{_RR}}$, the sequence (1) splits, which implies that $M\in
\mathcal{T}(\mod R)$ by Lemma 2.3. \hfill$\square$

\vspace{0.2cm}

Similarly, we have the following result.

\vspace{0.2cm}

{\bf Proposition 4.2.} {\it The following statements
are equivalent.

(1) ${^{\bot}{_RR}}\subseteq \mathcal{T}_1(\mod R)$. In this case,
${^{\bot}{_RR}}$ is said to have the torsionless property.

(2) ${^{\bot}{_RR}}\subseteq \mathcal{T}(\mod R)$.

(3) Every module in ${^{\bot}{_RR}}$ has finite torsionfree
dimension.

(4) $\mathcal{T}(\mod R^{op})\subseteq {^{\bot}R_R}$.}

\vspace{0.2cm}

Let $N \in\mod R^{op}$ and
$$0 \to N \buildrel {\delta _0} \over \longrightarrow E_{0}
\buildrel {\delta _1} \over \longrightarrow E_{1} \buildrel {\delta
_2} \over \longrightarrow \cdots \buildrel {\delta _i} \over
\longrightarrow E_{i} \buildrel {\delta _{i+1}} \over
\longrightarrow \cdots$$ be an injective resolution of $N$. For a
positive integer $n$, recall from [3] that an injective resolution
as above is called {\it ultimately closed} at $n$ if $\Im \delta
_n=\bigoplus _{j=0}^mW_j$, where each $W_j$ is a direct summand of
$\Im \delta _{i_j}$ with $i_j <n$. By [8, Corollary 2.3], if $R_R$
has a ultimately closed injective resolution at $n$ or $\id
_{R^{op}}R\leq n$, then ${^{\bot_n}{_RR}}$ (and hence
${^{\bot}{_RR}}$) has the torsionless property.

The following result generalizes [16, Lemma A], which states that
$\id _{R^{op}}R=\id _RR$ if both of them are finite.

\vspace{0.2cm}

{\bf Corollary 4.3.} {\it If $n=\min\{t\ |\ {^{\bot_t}{_RR}}$ has the
torsionless property$\}$ and $m=\min\{s\ |\ {^{\bot_s}{R_R}}$ has
the torsionless property$\}$, then $n=m$.}

\vspace{0.2cm}

{\it Proof.} We may assume that $n \leq m$. Let $N\in
{^{\bot_n}R_R}$. Then $N\in {^{\bot}R_R}(\subseteq {^{\bot_m}R_R})$
by Proposition 4.1. So $N\in \mathcal{T}(\mod R^{op})$ and
${^{\bot_n}R_R}$ has the torsionless property by the symmetric
version of Proposition 4.1. Thus $n \geq m$ by the minimality of
$m$. The proof is finished. \hfill$\square$

\vspace{0.2cm}

In the following, we will investigate the relation between the
torsionless property and the self-injective dimension of $R$. We
have seen that if $\id _{R^{op}}R\leq n$, then ${^{\bot_n}{_RR}}$
has the torsionless property. In the rest of this section, we will
investigate when the converse of this assertion holds true.







\vspace{0.2cm}

{\bf Proposition 4.4.} {\it Assume that $m$ and $n$ be positive
integers and $\Omega ^m (\mod R^{op})\subseteq$\linebreak
$\mathcal{T}_n (\mod R^{op})$. If $^{\bot _n}{_RR}$ has the
torsionless property, then $\id _{R^{op}}R\leq m$.}

\vspace{0.2cm}

{\it Proof.} Let $M\in \Omega ^m (\mod R^{op})$. Then $M\in
\mathcal{T}_n (\mod R^{op})$ by assumption. Because $^{\bot
_n}{_RR}$ has the torsionless property by assumption, $M\in
{^{\bot}{R_R}}$ by Proposition 4.1. Then it is easy to verify that
$\id _{R^{op}}R\leq m$. \hfill$\square$

\vspace{0.2cm}

Assume that
$$0\to {_RR}\to I^0(R)\to I^1(R) \to \cdots\to I^i(R)\to \cdots$$
is a minimal injective resolution of $_RR$.

\vspace{0.2cm}

{\bf Lemma 4.5.} {\it If $^{\bot _n}{_RR}$ has the torsionless
property, $\bigoplus_{i=0}^n I^i(R)$ is an injective cogenerator for
$\Mod R$.}

\vspace{0.2cm}

{\it Proof.} For any $S\in \mod R$, we claim that $\Hom _{R}(S,
\bigoplus_{i=0}^n I^i(R))\neq 0$. Otherwise, we have that $\Ext ^i
_{R}(S, R)\cong \Hom _{R}(S, I^i(R))=0$ for any $0\leq i \leq n$. So
$S\in {^{\bot _n}{_RR}}$ and hence $S$ is reflexive by assumption
and Proposition 4.1, which yields that $S\cong S^{**}=0$. This is a
contradiction. Thus we conclude that $\bigoplus_{i=0}^n I^i(R)$ is
an injective cogenerator for $\Mod R$. \hfill{$\square$}

\vspace{0.2cm}

{\bf Proposition 4.6.} {\it $\id _{R^{op}}R<\infty$ if and only if
$^{\bot _n}{_RR}$ has the torsionless property for some $n\geq 1$
and $\fd _R\bigoplus _{i\geq 0}I^i(R)<\infty$.}

\vspace{0.2cm}

{\it Proof.} The sufficiency follows from Lemmas 4.5 and 3.7, and
the necessity follows from Proposition 4.1 and Lemma 3.7.
\hfill{$\square$}

\vspace{0.2cm}

For any $n, k \geq 1$, recall from [9] that $R$ is said to be
$g_n(k)$ if $\Ext_{R^{op}}^j(\Ext_{R}^{i+k}(M, R), R))$ \linebreak
$=0$ for any $M\in \mod R$ and $1\leq i\leq n$ and $0\leq j\leq
i-1$; and $R$ is said to be $g_n(k)^{op}$ if $R^{op}$ is $g_n(k)$.
It follows from [12, 6.1] that $R$ is $g_n(k)$ (resp. $g_n(k)^{op}$)
if $\fd _{R^{op}}I^{i}(R^{op})$ (resp. $\fd _{R}I^{i}(R))\leq i+k$
for any $0\leq i\leq n-1$.

\vspace{0.2cm}

{\bf Theorem 4.7.} {\it Assume that $n$ and $k$ are positive integers
and $^{\bot_n}{_RR}$ has the torsionless property. If $R$ is
$g_n(k)$ or $g_n(k)^{op}$, then $\id _{R^{op}} R \leq n+k-1$.}

\vspace{0.2cm}

{\it Proof.} Assume that $^{\bot_n}{_RR}$ has the torsionless
property.

If $R$ is $g_n(k)$, then $\Omega ^{n+k-1}(\mod R)\subseteq
\mathcal{T}_{n}(\mod R)=\mathcal{T}(\mod R)$ by [9, Theorem 3.4]
and Proposition 4.1, which implies that the torsionfree dimension of
every module in $\mod R$ is at most $n+k-1$. So $\id _{R^{op}} R\leq
n+k-1$ by Theorem 3.6.

If $R$ is $g_n(k)^{op}$, then $\Omega ^{n+k-1} (\mod
R^{op})\subseteq \mathcal{T}_{n}(\mod R^{op})$ by the symmetric
version of [9, Theorem 3.4], which implies $\id _{R^{op}} R\leq
n+k-1$ by Proposition 4.4. \hfill$\square$

\vspace{0.2cm}

By Proposition 4.1 and Proposition 4.6 or Theorem 4.7, we
immediately get the following

\vspace{0.2cm}

{\bf Corollary 4.8.} {\it If $\fd _R\bigoplus _{i=0}^nI^i(R)\leq n$,
then $\id _{R^{op}} R \leq n$ if and only if $^{\bot _n}{_RR}$ has
the torsionless property.}

\vspace{0.2cm}

Recall that the Gorenstein symmetric conjecture states that $\id _R
R=\id _{R^{op}}R$ for any Artinian algebra $R$, which remains still
open. Hoshino proved in [5, Proposition 2.2] that if $\id _R R\leq
2$, then $\id _RR=\id _{R^{op}} R\leq 2$ if and only if $^{\bot
_2}{_RR}$ has the torsionless property. As an immediate consequence
of Theorem 4.7, the following corollary generalizes this result.

\vspace{0.2cm}

{\bf Corollary 4.9.} {\it For any $n\geq 1$, if $\id _R R\leq n$,
then $\id _R R=\id _{R^{op}}R\leq n$ if and only if $^{\bot
_n}{_RR}$ has the torsionless property.}

\vspace{0.2cm}

{\it Proof.} The necessity follows from Proposition 4.1. We next
prove the sufficiency. If $\id _R R\leq n$, then $\fd
_{R^{op}}\bigoplus _{i=0}^nI^{i}(R^{op})\leq n$ by Lemma 3.7, which
implies that $R$ is $g_n(n)$ by [12, 6.1]. Thus $\id _{R^{op}}\leq
2n-1$ by Theorem 4.7. It follows from [16, Lemma A] that $\id
_{R^{op}}R\leq n$. \hfill$\square$

\vspace{0.5cm}

\centerline{\bf 5. Questions}

\vspace{0.2cm}

In view of the results obtained above, the following two questions
are worth being studied.

Note that both the subcategory of $\mod R$ consisting of modules
with Gorenstein dimension at most $n$ and that consisting of modules
with left dimension at most $n$ are closed under extensions and
under kernels of epimorphisms. So, it is natural to ask the
following

\vspace{0.2cm}

{\bf Question 5.1.} {\it Is the subcategory of $\mod R$ consisting of
modules with torsionfree dimension at most $n$ closed under
extensions or under kernels of epimorphisms? In particular, Is
$\mathcal{T}(\mod R)$ closed under extensions or under kernels of
epimorphisms?}

\vspace{0.2cm}

For any $n\geq 1$, $\mathcal{T}_n(\mod R)$ is not closed under
extensions by [6, Theorem 3.3]. On the other hand, we have the
following

\vspace{0.2cm}

{\bf Claim.} {\it If ${^{\bot}{R_R}}$ has the torsionless property,
then the answer to Question 5.1 is positive.}

\vspace{0.2cm}

In fact, if ${^{\bot}{R_R}}$ has the torsionless property, then, by
the symmetric version of Proposition 4.2, we have that
$\mathcal{T}(\mod R)\subseteq {^{\bot}{_RR}}$ and every module in
$\mathcal{T}(\mod R)$ has Gorenstein dimension zero. So the
torsionfree dimension and the Gorenstein dimension of any module in
$\mod R$ coincide, and the claim follows.

\vspace{0.2cm}

By the symmetric version of [8, Corollary 2.3], if $_RR$ has a
ultimately closed injective resolution at $n$ or $\id _{R}R\leq n$,
then the condition in the above claim is satisfied. This fact also
means that the above claim extends [6, Corollary 2.5].

It is also interesting to know whether the converse of Theorem 3.6
holds true. That is, we have the following

\vspace{0.2cm}

{\bf Question 5.2.} {\it Does $\id _{R^{op}} R\leq n$ imply that
every module $M\in \mod R$ has torsionfree dimension at most $n$?}

\vspace{0.2cm}

{\bf Claim.} When $n=1$, the answer to Question 5.2 is positive.

\vspace{0.2cm}

Assume that $\id _{R^{op}} R\leq 1$ and $0\to K \to P \to M \to 0$
is an exact sequence in $\mod R$ with $P$ projective. Then $\Ext
_{R^{op}}^{i}(\Tr K, R)=0$ for any $i\geq 2$. Notice that $K$ is
torsionless, so $\Ext _{R^{op}}^{1}(\Tr K, R)=0$ and $K\in
\mathcal{T}(\mod R)$, which implies $\mathcal{T}-\dim _RM \leq 1$.
Consequently the claim is proved.

\vspace{0.5cm}

{\bf Acknowledgements.}
The authors thank the referee for the useful suggestions.

\newpage

\noindent Chonghui Huang\\
Research Institute of Mathematics\\
University of South China\\
Hengyang 421001, Hunan Province\\
P.R. China\\
and\\
Department of Mathematics\\
Nanjing University\\
Nanjing 210093, Jiangsu Province\\
P.R. China\\
e-mail: huangchonghui@usc.edu.cn

\vspace{0.5cm}

\noindent Zhaoyong Huang\\
Department of Mathematics\\
Nanjing University\\
Nanjing 210093, Jiangsu Province\\
P.R. China\\
e-mail: huangzy@nju.edu.cn

\end{document}